\documentclass[a4paper,english]{article}

\usepackage[english]{babel}
\usepackage[utf8]{inputenc}
\usepackage{hyperref}
\usepackage{amsmath}
\usepackage{amsfonts}
\usepackage{amssymb}
\usepackage{mathtools}
\usepackage{amsthm}
\usepackage{graphicx}
\usepackage{tikz}

\usetikzlibrary{patterns, shapes, arrows, positioning, fit, backgrounds, calc}

\usepackage{zibtitlepage}
\usepackage{authblk}

\definecolor{myblue}{RGB}{7, 127, 247}
\definecolor{myred}{RGB}{203, 65, 84}
\definecolor{mygreen}{RGB}{85, 107, 47}

\newtheorem{definition}{Definition}

\newtheorem{lemma}{Lemma}

\usepackage[style=apa,backend=biber]{biblatex}
\usepackage[a4paper, total={6in, 9in}]{geometry}

\begin{document}

\title{Enhancing Multi-Energy Modeling:\\ The Role of Mixed-Integer Optimization Decisions}

\author[1]{Stephanie Riedmüller\footnote{corresponding author, riedmueller@zib.de, \ZTPOrcid{0009-0006-4508-4262}}}
\author[1]{Annika Buchholz\footnote{buchholz@zib.de, \ZTPOrcid{0009-0000-2110-6929}}}
\author[1]{Janina Zittel\footnote{zittel@zib.de, \ZTPOrcid{0000-0002-0731-0314}}}

\affil[1]{Zuse Institute Berlin, Applied Algorithmic Intelligence Methods, Berlin, Germany}

\maketitle

\begin{abstract}

The goal to decarbonize the energy sector has led to increased research in modeling and optimizing multi-energy systems. One of the most promising and popular techniques for modeling and solving (multi-)energy optimization problems is (multi-objective) mixed-integer programming, valued for its ability to represent the complexities of integrated energy systems.
While the literature often focuses on deriving mathematical formulations and parameter settings, less attention is given to critical post-formulation decisions.

Modeling multi-energy systems as mixed-integer linear optimization programs demands decisions across multiple degrees of freedom. Key steps include reducing a real-world multi-energy network into an abstract topology, defining variables, formulating the relevant (in-)equalities to represent technical requirements, setting (multiple) objectives, and integrating these elements into a mixed-integer program (MIP). However, with these elements fixed, the specific transformation of the abstract topology into a graph structure and the construction of the MIP remain non-uniquely. These choices can significantly impact user-friendliness, problem size, and computational efficiency, thus affecting the feasibility and efficiency of modeling efforts.

In this work, we identify and analyze the additional degrees of freedom and describe two distinct approaches to address them. The approaches are compared regarding mathematical equivalence, suitability for solution algorithms, and clarity of the underlying topology. A case study on a realistic subarea of Berlin's district heating network involving tri-objective optimization for a unit commitment problem demonstrates the practical significance of these decisions. 

By highlighting these critical yet often overlooked aspects, our work equips energy system modelers with insights to improve computational efficiency, scalability, and interpretability in their optimization efforts, ultimately enhancing the practicality and effectiveness of multi-energy system models.
\end{abstract}

\section{Introduction}
Modeling and optimizing multi-energy systems has become a well-used approach to support the transformation of the energy infrastructure. 
Multi-energy systems combine multiple energy infrastructures (power, heat, cooling, fuel, transport) and technology regarding the production, consumption, storage, and transportation of energy (e.g., combined heat and power plants, heat pumps, steam turbines) into a holistic model, to gain insights about the links between several sectors. 

The literature on modeling multi-energy systems is rich; we refer to \textcite{Manco2024} for a general categorization.
Here, we are especially interested in multi-energy system optimization models implemented as a framework. In recent years, several approaches to categorizing and comparing the existing tooling landscape have been published, all of them considering different discriminating elements \parencite{Klemm2021,Kriechbaum2018,Hoffmann2024,Sola2020,Mancarella2014}.

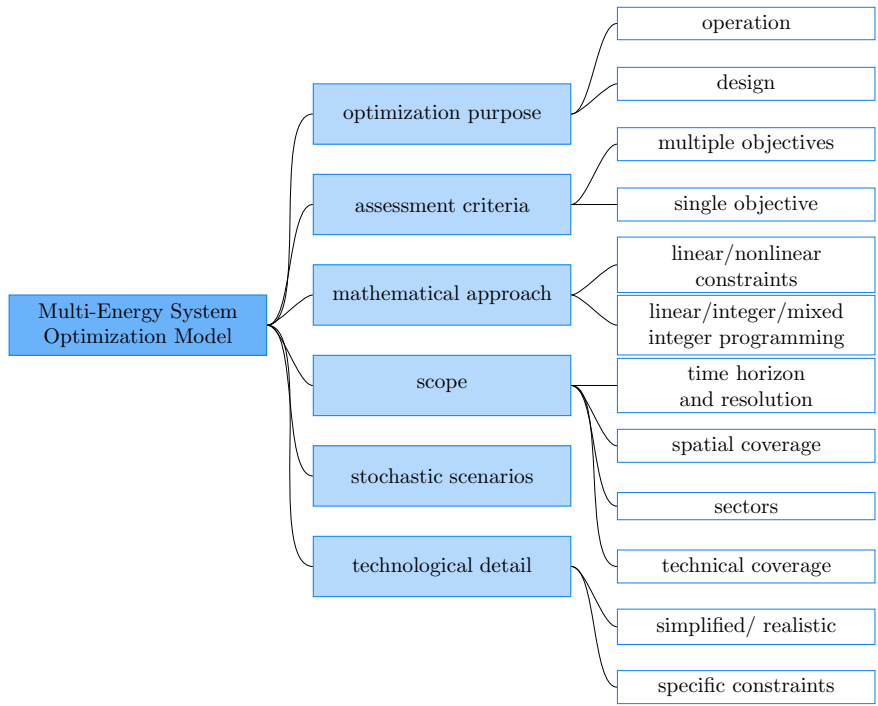
\begin{figure}[h!]
    \centering
    \scalebox{0.8}{\tikzstyle{topic} = [rectangle, 
minimum width=4cm, 
minimum height=1cm,
text centered, 
text width=4cm,
draw=myblue, 
fill=myblue!60]

\tikzstyle{category} = [rectangle, 
minimum width=4cm, 
minimum height=1cm, 
text centered, 
text width=4cm,
draw=myblue, fill=myblue!30]

\tikzstyle{subcategory} = [rectangle,  
minimum width=4cm, 
text centered, 
text width=4cm, 
draw=myblue, 
fill=white]
\tikzstyle{arrow} = [thick,->,>=stealth]

    \begin{tikzpicture}[node distance=1.5 cm]
    \node (mes) [topic] at (0.0, -0.5) {Multi-Energy System Optimization Model};
    
    \node (math) [category] at (5.0, 0.0) {mathematical approach};
    \node (lin) [subcategory] at (10.0, 0.5) {linear/nonlinear constraints};
    \node (pgrm) [subcategory] at (10.0, -0.5) {linear/integer/mixed integer programming};
    
    \node (criteria) [category,  above of=math] {assessment criteria};
    \node (single) [subcategory] at (10.0, 1.5) {single objective};
    \node (multiple) [subcategory] at (10.0, 2.5) {multiple objectives};
    
    \node (purpose) [category,  above of=criteria] {optimization purpose};
    \node (design) [subcategory] at (10.0, 3.5) {design};
    \node (operation) [subcategory] at (10.0, 4.5) {operation};

    \node (scope) [category,  below of=math] {scope};
    \node (space) [subcategory] at (10.0, -2.5) {spatial coverage};
    \node (time) [subcategory] at (10.0, -1.5) {time horizon and resolution};
    \node (sectors) [subcategory] at (10.0, -3.5) {sectors};
    \node (comp) [subcategory] at (10.0, -4.5) {technical coverage};

    \node (stoch) [category, below of=scope] {stochastic scenarios};
    
    \node (tech) [category,  below of=stoch]  {technological detail};
    \node (simpl) [subcategory] at (10.0, -5.5) {simplified/ realistic};
    \node (constr) [subcategory] at (10.0, -6.5) {specific constraints};

    \draw (mes) to [out=0, in=180, looseness=0.5] (math);
    \draw (math) to [out=0, in=180, looseness=0.5] (lin);
    \draw (math) to [out=0, in=180, looseness=0.5] (pgrm);
    \draw (mes) to [out=0, in=180, looseness=0.5] (criteria);
    \draw (criteria) to [out=0, in=180, looseness=0.5] (single);
    \draw (criteria) to [out=0, in=180, looseness=0.5] (multiple);
    \draw (mes) to [out=0, in=180, looseness=0.5] (purpose);
    \draw (purpose) to [out=0, in=180, looseness=0.5] (design);
    \draw (purpose) to [out=0, in=180, looseness=0.5] (operation);
    \draw (mes) to [out=0, in=180, looseness=0.5] (scope);
    \draw (scope) to [out=0, in=180, looseness=0.5] (space);
    \draw (scope) to [out=0, in=180, looseness=0.5] (time);
    \draw (scope) to [out=0, in=180, looseness=0.5] (sectors);
    \draw (scope) to [out=0, in=180, looseness=0.5] (comp);
    \draw (mes) to [out=0, in=180, looseness=0.5] (stoch);
    \draw (mes) to [out=0, in=180, looseness=0.5] (tech);
    \draw (tech) to [out=0, in=180, looseness=0.5] (simpl);
    \draw (tech) to [out=0, in=180, looseness=0.5] (constr);

    \end{tikzpicture}}
    \caption{Overview of the common discriminating elements for multi-energy system models and frameworks.}
    \label{fig:modelling_categories}  
\end{figure}

The following list represents some of the well-used dimensions for categorizing multi-energy systems optimization models and existing frameworks:
The \emph{optimization purpose} describes the main problem class, which usually involves a design and/or operation optimization. Defining the number and meaning of the objectives is part of the \emph{assessment criteria}. The \emph{mathematical approach} involves not only the general strategy, such as linear, nonlinear, or mixed-integer programming but also the (non)linearity of constraints and objectives. The \emph{scope} is defined by the coverage and resolution of the model regarding the time and spatial dimensions, considered energy sectors (heating sector, power sector, ...), and technical components.
We further distinguish between \emph{deterministic} and \emph{stochastic scenarios}. \emph{Technical details} can be on a range between realistic and simplified.
An overview is given in Figure~\ref{fig:modelling_categories}.

In this paper, the optimization technique is fixed to be a MIP for a unit commitment model, and the multi-energy system is based on a graph structure. In a previous study, we established a fixed mathematical formulation of the MIP with multiple objectives, which we evaluated to be suitable for a multi-energy system of large urban scale \parencite{Riedmueller2024, Zittel2024}. 

Even though the model is defined in the literature's categories, several dimensions in the implementation process are not uniquely determined and are usually not discussed. Those include decisions regarding the underlying graph structure and the transformation into the defined MIP. 
We note that the open dimensions considered are valid for other graph-based multi-energy system optimization problems. 
However, the computational study is restricted to the described case.

The goal is to compare the influence of specific decisions in the graph representation and their transfer into the corresponding MIP structure. The paper is structured as follows.
In Section~\ref{section:problem}, we describe the problem class, define the assumptions on which the investigation of this paper is based, and identify possible decision dimensions. 
Section~\ref{section:models} defines two model types, including amplitude decisions regarding the previously defined dimensions. Those two models stand diametrical regarding human readability and mathematical flexibility.
In Section~\ref{section:comparison}, we compare the influence of the decisions in the specific sub-categories regarding human readability, flexibility for the application of mathematical algorithms, and size of resulting MIP in theory.
Section~\ref{section:study} is dedicated to a computational investigation in a study on realistic instances of a subnetwork of Berlin's district heating network for different heating periods.

\section{Problem description and modeling base} \label{section:problem}

As a basis, we define the multi-energy system as consisting of a set of resources $R$, generating units $I$, and storage units $K$. Further, we define a set of markets $M$ (e.g. fuel and power markets) , demand nodes $D$ and balance nodes $B$. We consider the network to be represented by a graph $G = (V, A)$ with 
$V = I \cup K \cup M \cup D \cup B$
and 
$A$ containing the arcs given by the topology of the network.
All inputs are given as a series of discretized time steps, where $T$ denotes the set of time steps.
For this graph structure, we want to solve operation optimization. Operation optimization can be modeled as a unit commitment problem, which involves
\begin{itemize}
    \item commitment decision (whether a unit is operated at a given time step),
    \item production decision (how much energy is produced by a unit at a given time step),
    \item network flow decision (how much energy is flowing on an arc at a given time step).
\end{itemize}

The variables $x_{t, v^{in}}^r, x_{t, v^{out}}^r \in \mathbb{R}_{\geq 0}$ denote the incoming and outgoing flow of resource $r \in R$ at node $v \in V$ and time $t \in T$. Note that for balance nodes the incoming and outgoing flow can be a sum of flows. The status of generating unit $i$ and time $t$ is denoted by $z_{i,t} \in \{0,1 \}$. Variable $s_{i,t} \in \{0,1\}$ defines, if a status change has been executed from time step $t-1$ to $t$. The storage level of storage $k \in K$ is denoted by $h^r_{t,k}$ and the variables $p_t^r, e_t^r \in \mathbb{R}_{\geq 0}$ represent the purchased and sold resources. 

We minimize operational costs, CO$_2$ emissions, and maximize the produced heat of combined heat and power plants (CHP). For given demand vector $d$, parameter vectors $a$ and $c$, and for a conversion map $\varphi$, the following MIP models operation optimization for multi-energy systems:
\newline
\resizebox{\linewidth}{!}{
\begin{minipage}{1.05\linewidth}
\begin{align}
    &\text{min} && \left(\text{costs, emissions, $-$ CHP heat} \right) &&& \label{mip}\tag{$MIP$}\\
    & && \sum_{i \in I}  x_{t, i^{out}}^{r} + \sum_{k \in K} x_{t, k^{out}}^{r} + p_t^r  = d^r + \sum_{i \in I}  x_{t, i^{in}}^{r} + \sum_{k \in K} x_{t, k^{in}}^{r} + e_t^r &&& \forall r \in R \label{eq:balance}\\
    & && x_{t, i^{out}}^{r_2} = \varphi_{i,t}^{r_1,r_2} \left(x_{t, i^{in}}^{r_1}\right) &&& \forall i \in I, t\in T, r_1, r_2 \in R \label{eq:conversion}\\
    & && s_{i,t} \leq z_{i,t}, s_{i,t} \leq 1 - z_{i,t}, s_{i,t} \geq z_{i,t} - z_{i,t-1} &&& \forall i \in I, t\in T \label{eq:activation}\\
    & && \sum_{\tau \in T^{up}_{i,t}} ( s_{i,t} - z_{i,\tau}) \leq 0, \sum_{\tau \in T^{down}_{i,t}} ( s_{i,t} + z_{i,\tau} -1) \leq 0 &&& \forall i \in I, t\in T \label{eq:minupdown}\\
    & && x_{t+1, i^{out}}^{r} - x_{t, i^{out}}^{r} \leq a_i^{up}, x_{t, i^{out}}^{r} - x_{t+1, i^{out}}^{r} \leq a_i^{down} &&&\forall i \in I, t\in T, r \in R \label{eq:ramping}\\
    & && h^r_{t+1,k} = a^{loss}_{t,k}h^r_{t,k} + a^{load}_{t,k}x_{t, k^{in}}^{r} - a^{unload}_{t,k}x_{t, k^{out}}^{r} &&&\forall k \in K, t\in T, r \in R \label{eq:storage}\\
    & && h,x,p,e \leq c_{max}, \qquad h,x,p,e \geq c_{min} \label{eq:capacities} \\
    \nonumber
\end{align}
\end{minipage}
}

Constraint (\ref{eq:balance}) represents demand and resource balance.
Constraint (\ref{eq:conversion}) defines the resource conversion at a node through a conversion map $\varphi$, which usually describes a technology-specific characteristic curve. Here, we define the $\varphi$ to be a piecewise linear function. 
Constraints (\ref{eq:activation}) - (\ref{eq:ramping}) includes technical constraints (activation, minimum up and down times, ramping), 
Constraint (\ref{eq:storage}) time linking of storage levels and
Constraint (\ref{eq:capacities}) capacities.
For a in-depth explanation of the model, see \textcite{Clarneretal2020}.

Even with the fixed choices, there are still decisions to be made that are given only implicitly by the described model. Those decisions are often not mentioned but can play a role in user-friendliness and computational efficiency.
We identify the following open dimensions:
On the one hand, explicit or implicit representation of components in the underlying graph structure can have a significant impact. The underlying graph structure could be solely based on the topology of the main assets of the considered real-world network instance or include components used for modeling purposes, even if they are not present as real-world assets. Those objects include the representation of a hierarchical level structure, the degree of specification of node types, the explicit representation of subcomponents, and the modeling of individual objective functions, for example, in the transportation of information.
On the other hand, it is a significant design decision if variables in the MIP are based on nodes or arcs in the underlying graph. Even though variables are already defined in the previous problem description, the realization of the stated MIP is not uniquely determined.

\section{Description of model types} \label{section:models}

We describe two diametrically engineered models regarding the previously stated dimensions in the following. Model A prioritizes human readability, while Model B prioritizes mathematical flexibility. Table~\ref{tab:overview_dimensions} depicts an overview of the characteristics. Both options are valid candidates as their characteristics are widely used in modeling communities. Exemplary frameworks are given for the respective characteristics of both models. We consider the openly available graph-based multi-energy system frameworks Calliope \parencite{calliope}, urbs \parencite{urbs}, oemof \parencite{oemof} and PyPSA \parencite{PyPSA} and the commercial software BoFiT \parencite{bofit}.

\begin{table}[h!]
    \centering
    \begin{tabular}{lllll}
    \hline
           & Dimension         & Model A  & Model B\\
    \hline
        Graph representation   & structure   & hierarchical & homogeneous\\
                    & information transport & explicit  & implicit  \\
                    & specification of node types        & technology specific & abstract\\
                    &  connections       & defined ports & flexible\\
                    & subcomponents        & detailed & contraction\\
    \hline
        MIP         & association of flow variables    & nodes & arcs\\
                    & constraints for specific technologies    & detailed & reduced\\
    \hline
    \end{tabular}
    \caption{Overview of the characteristics of the various dimensions for Model A and Model B.}
    \label{tab:overview_dimensions}
\end{table}

\subsection{Model A: Prioritizing human readability}

Model A focuses on a user-friendly structure, which could, for example, be used to implement an energy modeling framework designed for commercial use.

\begin{figure}[h!]
    \centering
    \begin{minipage}{0.30\textwidth}
        \centering
        \scalebox{0.8}{\tikzstyle{container} = [rectangle, 
minimum width=1.5cm, 
minimum height=0.8cm, 
text centered, 
draw=myblue, fill=myblue!30]
\tikzstyle{container2} = [rectangle, 
minimum width=1.5cm, 
minimum height=0.8cm, 
text centered, 
draw=myblue]
\tikzstyle{containerlarge} = [rectangle, 
minimum width=10.7cm, 
minimum height=5.3cm, 
text centered, 
draw=myblue]
\tikzstyle{oval} = [rectangle, 
minimum width=1.5cm, 
minimum height=0.8cm, 
text centered, 
draw=myblue,
rounded corners=5pt]
\tikzstyle{enhanced} = [circle, 
minimum width=1cm, 
minimum height=1cm, 
text centered, 
dashed,
draw=myblue]
\tikzstyle{dummy} = [circle, 
text centered]
\tikzstyle{arrow} = [thick,->,>=stealth]

    \begin{tikzpicture}[node distance=2.3 cm]
    
    \node (fuel) [container] {Fuel Market};
    \node (site1) [container2,  above right of=fuel] {Site 1};
    \node (site2) [container,  below right of=fuel] {Site 2};
    \node (heat) [container, right of=site1] {Heat Export};
    \node (power) [container, right of=site2] {Power Export};
    
    \draw [arrow, myblue] (fuel) --  (site1);
    \draw [arrow, myblue] (site1) --  (heat);
    \draw [arrow, myblue] (site1) --  (power);
    \draw [arrow, myblue] (fuel) --  (site2);
    \draw [arrow, myblue] (site2) --  (power);

    \end{tikzpicture}}
    \end{minipage}
    \hfill
    \begin{minipage}{0.69\textwidth}
        \centering
        \scalebox{0.7}{\tikzstyle{container} = [rectangle, 
minimum width=1.5cm, 
minimum height=0.8cm, 
text centered, 
draw=myblue, fill=myblue!30]
\tikzstyle{containerlarge} = [rectangle, 
minimum width=9.2cm, 
minimum height=4.8cm, 
text centered, 
draw=myblue]
\tikzstyle{oval} = [rectangle, 
minimum width=1.0cm, 
minimum height=0.6cm, 
text centered, 
draw=myblue,
rounded corners=5pt]
\tikzstyle{enhanced} = [circle, 
minimum width=1cm, 
minimum height=1cm, 
text centered, 
dashed,
draw=myblue]
\tikzstyle{dummy} = [circle, 
text centered]
\tikzstyle{arrow} = [thick,->,>=stealth]

    \begin{tikzpicture}[node distance=2.2 cm]
        
    \node (in) [dummy] {};
    \node (balance1) [dummy, right of=in] {};
    \node (gas1) [oval, above right of=balance1] {\textcolor{myblue}{Gas Turbine 1}};
    \node (gas2) [oval, below right of=balance1] {\textcolor{myblue}{Gas Turbine 2}};
    \node (boiler) [oval, above right of=gas2] {\textcolor{myblue}{Heat Boiler}};
    \node (steam) [oval, above right of=boiler] {\textcolor{myblue}{Steam Turbine}};
    \node (converter) [oval, right of=steam] {\textcolor{myblue}{Converter}};
    \node (balance2) [dummy, right of=converter] {};
    \node (balance3) [dummy, below of=balance2] {};
    \node (out1) [dummy, right of=balance2] {};
    \node (out2) [dummy, right of=balance3] {};
    \node (boxdummy) [dummy, right of=balance1] {};
    \node (box) [containerlarge, right of=boxdummy] {};
    \node (namedummy) [dummy, right of=gas2] {};
    \node (name) [dummy, right of=namedummy] {Site 1};

    \draw [arrow, myblue, dashed] (in) -- (balance1);
    \draw [arrow, myblue] (balance1) -- (gas1);
    \draw [arrow, myblue] (balance1) -- (gas2);
    \draw [arrow, myblue] (gas1) -- (boiler);
    \draw [arrow, myblue] (gas2) -- (boiler);
    \draw [arrow, myblue] (boiler) -- (steam);
    \draw [arrow, myblue] (steam) -- (converter);
    \draw [arrow, myblue] (converter) -- (balance2);
    \draw [arrow, myblue] (converter) -- (balance3);
    \draw [arrow, myblue] (boiler) -- (balance3);
    \draw [arrow, myblue, dashed] (balance2) -- (out1);
    \draw [arrow, myblue, dashed] (balance3) -- (out2);

    \end{tikzpicture}}
    \end{minipage}
    \caption{Hierarchical structure in form of a visually contractable subgraph (Site 1).}\label{fig:hierarchical_structure}
\end{figure}
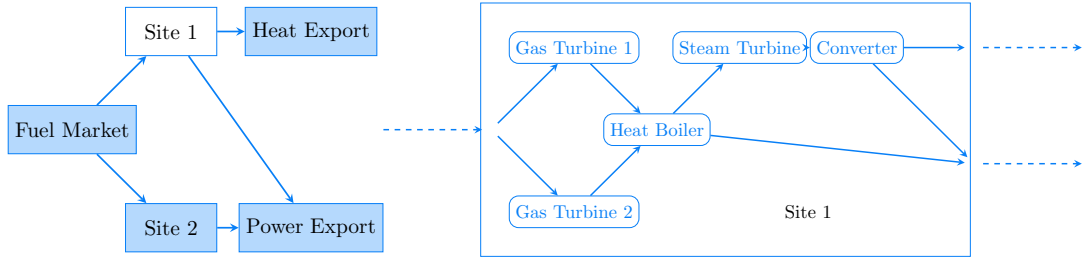

\paragraph{Hierarchical structure in graph representation:} 
 The structure of Model A allows the inclusion of hierarchical graph levels in the form of containers. 
A container node can encapsulate a subgraph. This allows a concise picture, especially if modeled in visual software, where users can reduce the visible network to specific needs.  
To enable this behavior, the input and output of such a container need to be defined, which connects the graph with the subgraph inside the container by exactly two arcs. Consider the example in Figure~\ref{fig:hierarchical_structure} depicting the hierarchical encapsulation applied to the different generating unit types combined in a site container.  A hierarchical structure in the form of energy hubs is, for example, used in Calliope and urbs. For the concept of energy hubs, we refer to \textcite{Kriechbaum2018}. PyPSA allows the explicit definition of subnetworks. A graphical interface to create models is given in BoFiT.

\begin{figure}[h!]
    \centering
    \scalebox{0.8}{\tikzstyle{containerfill} = [rectangle, 
minimum width=1.5cm, 
minimum height=0.8cm, 
text centered, 
draw=myblue, fill=myblue!30]
\tikzstyle{container} = [rectangle, 
minimum width=1.5cm, 
minimum height=0.8cm, 
text centered, 
draw=gray]
\tikzstyle{containerlarge} = [rectangle, 
minimum width=10.7cm, 
minimum height=5.3cm, 
text centered, 
draw=myblue]
\tikzstyle{oval} = [rectangle, 
minimum width=1.5cm, 
minimum height=0.8cm, 
text centered, 
draw=myblue,
rounded corners=5pt]
\tikzstyle{enhanced} = [circle, 
minimum width=1cm, 
minimum height=1cm, 
text centered, 
dashed,
draw=myblue]
\tikzstyle{dummy} = [circle, 
text centered]
\tikzstyle{arrow} = [thick,->,>=stealth]

    \begin{tikzpicture}[node distance=2.3 cm]
    
    \node (fuel) [container] {Fuel Market};
    \node (chp) [container,  right of=fuel] {CHP};
    \node (heat) [container,  above=1cm of chp] {Heat Demand};
    \node (power) [container, right of=chp] {Power Market};
    \node (cost) [containerfill, below=1cm of chp] {Costs};
    \node (dummy) [dummy, right of=cost] {};
    \node (objective) [containerfill, right of=dummy] {Objective Function};
    
    \draw [arrow] (fuel) --  (chp);
    \draw [arrow] (chp) --  (heat);
    \draw [arrow] (chp) --  (power);
    \draw [arrow, myblue] (cost) --  (chp);
    \draw [arrow, myblue] (fuel) --  (cost);
    \draw [arrow, myblue] (cost) --  (power);
    \draw [arrow, myblue] (power) --  (cost);
    \draw [arrow, myblue] (chp) --  (cost);
    \draw [arrow, myblue] (cost) --  (objective);

    \end{tikzpicture}}
    \caption{Cost-related information for the objective function is explicitly transported through the graph.}
    \label{fig:explicit_information}  
\end{figure}
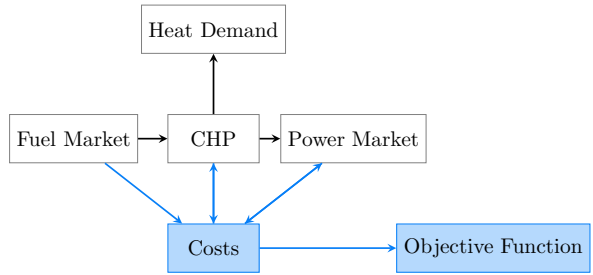

\paragraph{Explicit representation of information transport:}
Model A does not differ between the transport of actual commodities (heat, power, \ldots) and the transport of information (costs, temperature, \ldots). Therefore, all transport is explicitly represented in the graph structure. This approach is often used in highly visual software so that all configurations are visually accessible by the user. 
For example, multiple objectives are handled as artificial nodes containing ports for each objective and its positive and negative values. The artificial objective nodes are connected through arcs and possible conversion nodes to those nodes that influence the corresponding objective. This way, the information for the objectives is explicitly transported through the network, similar to other commodities, see Figure~\ref{fig:explicit_information}. On the one hand, this allows for easy handling of an objective, containing only variables given as ports of the objective node without additional scalarization since scalarization is already implicitly handled. On the other hand, the components of an objective are not easily accessible since they are already aggregated throughout the network before reaching the objective node. Information transport is represented explicitly, for example, in BoFiT.

\begin{figure}[h!]
    \centering
    \scalebox{0.8}{\tikzstyle{container} = [rectangle, 
minimum width=1.5cm, 
minimum height=0.6cm, 
text centered, 
draw=myblue,
fill=myblue!30]
\tikzstyle{largecontainer} = [rectangle, 
minimum width=10cm, 
minimum height=4.9cm, 
text centered, 
draw=myblue]
\tikzstyle{oval} = [rectangle, 
minimum width=1.5cm, 
minimum height=0.8cm, 
text centered, 
draw=myblue,
rounded corners=5pt]
\tikzstyle{enhanced} = [circle, 
minimum width=1cm, 
minimum height=1cm, 
text centered, 
dashed,
draw=myblue]
\tikzstyle{dummy} = [circle, 
text centered]
\tikzstyle{dummy2} = [rectangle, 
text centered]
\tikzstyle{arrow} = [thick,-,>=stealth]

    \begin{tikzpicture}[node distance=1.0 cm]

    \node (input) [container] {Actual Input};
    \node (d1) [dummy, left=0.5cm of input] {};
    \node (d2) [dummy, right=0.5cm of input] {};
    
    \node (Tin) [container, below of=input] {Thermal Input};
    \node (d3) [dummy, left=0.5cm of Tin] {};
    \node (d4) [dummy, right=0.5cm of Tin] {};
    
    \node (flow) [container, below of=Tin] {Flow Temperature};
    \node (d5) [dummy, left=0.5cm of flow] {};
    \node (d6) [dummy, right=0.5cm of flow] {};
    
    \node (opload) [container, below of=flow] {Operation Loading};
    \node (d7) [dummy, left=0.5cm of opload] {};
    \node (d8) [dummy, right=0.5cm of opload] {};
    
    \node (load) [container, below of=opload] {Loading};
    \node (d9) [dummy, left=0.5cm of load] {};
    \node (d10) [dummy, right=0.5cm of load] {};
    
    \node (dummy) [dummy, right=4cm of input] {};
    \node (storage) [dummy2, above of=dummy] {Storage};

    \node (output) [container, right=8cm of input] {Actual Output};
    \node (d11) [dummy, left=0.5cm of output] {};
    \node (d12) [dummy, right=0.5cm of output] {};
    
    \node (Tout) [container, below of=output] {Thermal Output};
    \node (d13) [dummy, left=0.5cm of Tout] {};
    \node (d14) [dummy, right=0.5cm of Tout] {};
    
    \node (return) [container, below of=Tout] {Return Temperature};
    \node (d15) [dummy, left=0.5cm of return] {};
    \node (d16) [dummy, right=0.5cm of return] {};
    
    \node (opunload) [container, below of=return] {Operation Unloading};
    \node (d17) [dummy, left=0.5cm of opunload] {};
    \node (d18) [dummy, right=0.5cm of opunload] {};
    
    \node (unload) [container, below of=opunload] {Unloading};
    \node (d19) [dummy, left=0.5cm of unload] {};
    \node (d20) [dummy, right=0.5cm of unload] {};

    \node (level) [container, below=3.9cm of dummy] {Storage Level};
    \node (d22) [dummy, above=0.5cm of level] {};
    \node (d23) [dummy, below=0.5cm of level] {};

    \begin{scope}[on background layer]
    \node (d21) [largecontainer, below=0.28cm of storage] {};
    \end{scope}
    
    \node (eq1) [dummy2, below=1.0cm of storage] {Temperature Constraints};
    \node (eq2) [dummy2, below=0.3cm of eq1] {Capacity Constraints};
    \node (eq3) [dummy2, below=0.3cm of eq2] {Storage Constraints};
;

    \draw [arrow, myblue] (d1) --  (input);
    \draw [arrow, myblue] (input) --  (d2);
    \draw [arrow, myblue] (d3) --  (Tin);
    \draw [arrow, myblue] (Tin) --  (d4);
    \draw [arrow, myblue] (d5) --  (flow);
    \draw [arrow, myblue] (flow) --  (d6);
    \draw [arrow, myblue] (d7) --  (opload);
    \draw [arrow, myblue] (opload) --  (d8);
    \draw [arrow, myblue] (d9) --  (load);
    \draw [arrow, myblue] (load) --  (d10);

    \draw [arrow, myblue] (d11) --  (output);
    \draw [arrow, myblue] (output) --  (d12);
    \draw [arrow, myblue] (d13) --  (Tout);
    \draw [arrow, myblue] (Tout) --  (d14);
    \draw [arrow, myblue] (d15) --  (return);
    \draw [arrow, myblue] (return) --  (d16);
    \draw [arrow, myblue] (d17) --  (opunload);
    \draw [arrow, myblue] (opunload) --  (d18);
    \draw [arrow, myblue] (d19) --  (unload);
    \draw [arrow, myblue] (unload) --  (d20);

    \draw [arrow, myblue] (d22) -- (level);
    \draw [arrow, myblue] (level) -- (d23);

    \end{tikzpicture}}
    \caption{Visualization of fixed ports belonging to a unit node on the example of a storage unit.}
    \label{fig:storage}  
\end{figure}
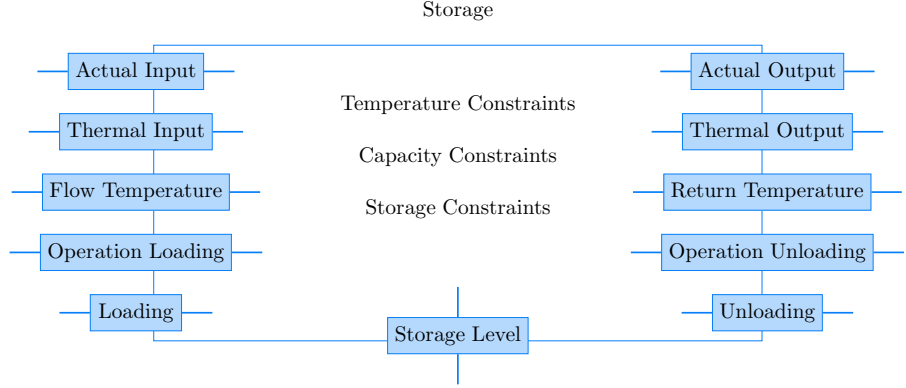

\paragraph{Detailed specification of nodes:}
In Model A, each specific technology is represented by a corresponding unit node. Each unit node includes technology-specific ports, i.e., specified incoming and outgoing arcs defined in recourse and number, see Figure~\ref{fig:storage}. If a generating unit consists of several subunits, all subunits are modeled individually with corresponding conversion information; see again Figure~\ref{fig:hierarchical_structure}. Furthermore, each technology has a corresponding set of constraints modeling the specific relations in detail.
Due to the exact modeling of the specific technologies, some constraints inside those specific technology unit nodes can be redundant. 
While containing rather abstract component types, PyPSA provides many detailed options for each component.
Components with fixed ports, for example, can be found in BoFiT. 

\begin{figure}[h!]
    \centering
    \tikzstyle{classic} = [circle, 
minimum width=0.5cm, 
minimum height=0.5cm, 
text centered, 
draw=myblue]
\tikzstyle{arrow} = [thick,->,>=stealth]

    \begin{tikzpicture}[node distance=4 cm]
    \node (n11) [classic] {v};
    \node (n22) [classic, right of=n11] {w};
    
    \node (n1) [classic] at (6,0) {v};
    \node (n2) [classic, right of=n1] {w};

    \draw [arrow] (n1) -- node [above,midway] {$x_{v,w}$} (n2);
    \draw [arrow] (n11) -- node [above,pos=.1] {$x_{v,out}$} node [above,pos=.9] {$x_{w,in}$} node [above,pos=.5] {$=$}(n22);

    \end{tikzpicture}
    \caption{Choice of flow variables in Model A (left) and Model B (right).}
    \label{fig:flow_vars}  
\end{figure}
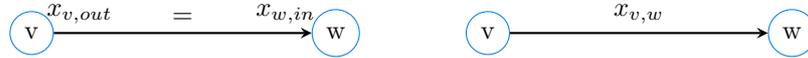

\paragraph{Focus on node functions:}
A major decision is how to model the flow variables $X$ given in the described MIP. 
Despite clearly modeling the flow of a specific resource on a given arc, the flow is often viewed as an incoming and outgoing resource of a node. Model A takes this viewpoint by relating the incoming and outgoing resources to the ports belonging to a node; see Figure~\ref{fig:flow_vars}. Ports can handle incoming and outgoing arcs and allow different possible requirements for arcs connected to the node. For example, a combined heat and power plant must have an incoming fuel arc and outgoing heat and power edges. If those ports are fixed, the function of different connected arcs is transparent, and connected arcs can be handled explicitly in the component. Hence, the flow on an arc $(v,w)$ is encoded in the incoming port variable of $w$ and the outgoing port variable of $v$. This approach allows for a node-focused model but needs an additional constraint to relate the variables on the ports belonging to the same arc. For example, Calliope and urbs model the nodes' energy flows as input and output.

Modal A describes an approach that results in a detailed depiction of network instances and technology often used for engineering purposes.

\subsection{Model B: Prioritizing mathematical flexibility}

Model B focuses on a restriction to a mathematically necessary base, allowing for flexibility for the application of further solving algorithms.

\paragraph{Homogeneous structure in graph representation:} 
Model B does not allow for a hierarchical level structure of the graph but is restricted to a universal base layer. 
This approach lacks a visually structured overview and is, therefore, without a further abstraction layer not suitable for visual software. But it leads to a clear, minimalistic, and homogeneous representation of the given information. In addition, the absence of additional layers reduces the graph size immensely and, therefore, the number of variables and flow conservation constraints (see Section~\ref{section:comparison_graph} and \ref{section:comparison_mip}). Frameworks usually offer hierarchical structures as an optional feature, allowing for a homogeneous structure if desired.

\paragraph{Implicit representation of information transport:}
In Model B, the transport of information that is not a commodity is only given implicitly, i.e., it has no explicit representation in the graph structure. The graph solely includes the topology of the network.
Multiple objectives are constructed in the MIP model outside of the graph structure by collecting relevant information from the graph after building the graph, including its MIP variables and constraints. This allows an exact understanding of the objective's construction and the corresponding scalarization of the variables. 
Due to this knowledge, algorithms handling the objectives, such as multi-objective algorithms as a weighted sum, can be applied in an informed manner. Furthermore, the results of optimization procedures can be interpreted more directly, and alterations of scalarizations can be applied easily. Implicit information transport can be found in the frameworks Calliope, urbs, and oemof.

\paragraph{Abstract nodes:}
The nodes of Model B are strongly abstracted in the form of general input/output components. Instead of specific numbers and resources in the incoming ports, the general input/output components are only specified in relation to the number of incoming and outgoing arcs. 
Technology-specific constraints are only included if necessary for the solution polytope, i.e., redundant constraints are excluded.

The conversion of resources in a node is given purely abstract as a piece-wise linear function.
If feasible, subcomponents are merged by integrating their inherent conversion processes into a single process. In addition to reducing the variables and constraints needed to model several conversion processes, we also save additional information transfer.
Some nodes convert, for example, the input fuel into heat or both heat and power. To reduce the number of nodes and constraints, any connected nodes that perform these functions can be merged into a single node. Merging can be beneficial in specific scenarios. For instance, if the first node is a fuel-to-heat converter and the second node is either a fuel-to-heat converter or a fuel-to-heat and power converter. Conversely, if the first node is a fuel-to-both converter, it can only be merged with a fuel-to-heat converter. The new converter must be defined using the incoming arcs of the first converter and the outgoing arcs of the second, reducing the complexity of Model B. When this merging occurs, the characteristic curve of the combined node must be adjusted accordingly:
The conversion curve is defined as a piece-wise linear function $\varphi$, represented by a source series $ x_0, x_1, \ldots, x_n $ and a target series $ y_0, y_1, \ldots, y_n $, where $y_i < y_j$ for $i<j$ holds. The condition for the curve to continuously increase is necessary for defining the decomposition of two or more curves. Let the second converter be given by the piecewise linear function $\psi$, defined by a source series $ \overline{y}_0, \overline{y}_1, \ldots, \overline{y}_n $ and a target series $ z_0, z_1, \ldots, z_n $, $z_i < z_j$ for $i < j$. The decomposition $ \varphi \circ \psi $ can be calculated by identifying the values $ z_i^* = \psi(y_i) $ and $ x_i^* = \varphi^{-1}(\overline{y}_i) $ for every $ i = 0, \ldots, n $. The resulting source and target series are given by the pairs $ (x_i, z^*_i) $ and $ (x^*_i, z_i) $, which must be sorted so that the source values are in increasing order.

The approach of preferably abstract nodes allows for a generalized depiction and handling of the MIP variables and leads to a more concise and abstract structure, as can be seen, for example, in urbs.

\paragraph{Focus on arc flows:}
Model B considers the flow variables X in a theoretical network flow manner, i.e., the flow is not viewed as a feature of the nodes but as a feature belonging to the arcs, see Figure~\ref{fig:flow_vars}. Hence, we generate exactly one flow variable for each arc in the graph. Since the nodes still include the conversion procedures of the resource flow, constraints modeling the conversion need access to the adjacent arc variables. This approach is used for modeling network flows in many graph theoretical algorithms. Those algorithms can, therefore, be applied with less effort to the model. The energy flow is modeled on arcs, for example, in oemof and PyPSA.

Model B is built to restrict the size of the graph and MIP to a minimum and allow for the application of mathematical algorithms.

\section{Comparison} \label{section:comparison}

In this section, we compare the two models across multiple, both qualitative and quantitative, dimensions to evaluate their structural and computational advantages. By soft criteria, we examine practical aspects such as interpretability and degree of detail. We further compare the graph size of the models and analyze the number of variables and constraints in the corresponding MIP for each model. Furthermore, we investigate mathematical equivalence, identifying decisions under which both models yield identical solutions. This comprehensive comparison provides insights into the trade-offs and suitability of each model for different applications.

\subsection{Soft criteria}

The hierarchical structure in the form of containers in Model A can increase the human readability of the network structure, especially for visual software. However, it is not necessary from the point of computation and can, in fact, complicate the application of algorithms. 

Modeling information transport explicitly as components in the underlying graph structure, as in Model A, is of interest if the instances are built in software using a visual representation. If the underlying graph structure does reflect solely an abstraction of the topology of the considered network, as in Model B, it offers a more intuitive and comprehensible overview.

The specification of components for different technologies, as given in Model A, gives the instances an application-near appearance but is not necessary for computations. However, an abstraction of similar components into a generic input/output component, as given in Model B, can be beneficial for formulating algorithms.
Including subcomponents explicitly in the graph structure, as in Model A, can give a more detailed representation of functionality and purpose comparable to the network's topology. However, combining several detailed subcomponents, as in Model B, can increase the overview. 

Furthermore, Model A results in a closure of individual components, which allows constraints to be formulated on a node completely independent of the incident arcs. In contrast, the modeling of individual components in Model B relies on the availability of information of the incident arcs. This difference must be considered when deciding for a data structure used in an implementation.

\subsection{Comparison of graph size} \label{section:comparison_graph}

\begin{figure}[h!]
    \centering
    \tikzstyle{container} = [rectangle, 
minimum width=5.2cm, 
minimum height=3cm, 
text centered, 
draw=myblue]
\tikzstyle{classic} = [circle, 
minimum width=0.5cm, 
minimum height=0.5cm, 
text centered, 
draw=black, fill=gray!30]
\tikzstyle{enhanced} = [circle, 
minimum width=0.5cm, 
minimum height=0.5cm, 
text centered, 
dashed,
draw=myblue]
\tikzstyle{dummy} = [circle, 
text centered]
\tikzstyle{arrow} = [thick,->,>=stealth]

    \begin{tikzpicture}[node distance=1.5 cm]
    \node (n0) [dummy] {};
    \node (n1) [enhanced,  right of=n0] {};
    \node (n2) [classic,  right of=n1] {};
    \node (n3) [classic,  above right of=n2] {};
    \node (n4) [classic,  below right of=n2] {};
    \node (n5) [classic,  right of=n3] {};
    \node (n6) [classic,  below right of=n5] {};
    \node (n7) [enhanced,  right of=n6] {};
    \node (n8) [dummy,  right of=n7] {};
    \node (cont) [container,  right= 2.0cm of n0] {};
    \node (name) [dummy,  below right of=cont] {\textcolor{myblue}{Container $c \in C$}};

    \draw [arrow] (n0) --  (n1);
    \draw [arrow, myblue, dashed] (n1) --  (n2);
    \draw [arrow] (n2) --  (n3);
    \draw [arrow] (n2) --  (n4);
    \draw [arrow] (n3) --  (n5);
    \draw [arrow] (n4) --  (n6);
    \draw [arrow] (n5) --  (n6);
    \draw [arrow] (n4) --  (n5);
    \draw [arrow, myblue, dashed] (n6) --  (n7);
    \draw [arrow] (n7) --  (n8);

    \node[] at (0,1) {
    \begin{tikzpicture}[node distance=2mm and 1mm]

      \node (dummy) [white, minimum width=4mm, minimum height=4mm] {};
      
      \node (sq1) [fill=black, minimum width=4mm, minimum height=4mm, below= of dummy] {};
      \node [right=of sq1] {Model B};

      \begin{scope}
        \coordinate (base) at (-0.2,-0.2);
        \fill[myblue] (base) rectangle ++(2mm,4mm);
        \fill[black]  ($(base)+(2mm,0mm)$) rectangle ++(2mm,4mm); 
        \node [right=of dummy] {Model A};
      \end{scope}
    \end{tikzpicture}
  };

    \end{tikzpicture}
    \caption{Depicted is a subgraph in Model A with a hierarchical level structure in form of a container (illustrated in blue), which allows the subgraph to be contracted in a visual software. The black graph shows the corresponding subgraph in Model B with a homogeneous structure. The dotted, blue parts represent additionally generated nodes and arcs in Model A in comparison to Model B. Compare Figure ~\ref{fig:hierarchical_structure}.}
    \label{fig:container}  
\end{figure}
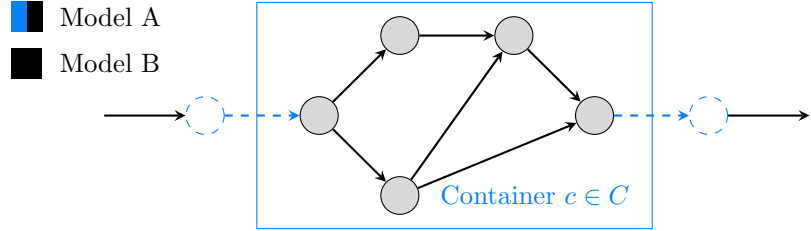

For each container structure, Model A includes two balance nodes and two arcs more than Model B, as depicted in Figure~\ref{fig:container}. Hence, in Model A, the number of nodes and arcs in the underlying graph is increased by $2 \cdot |C|$ compared to Model B, where $C$ denotes the set of all container structures $c$ in Model A.

Regarding the explicit or implicit representation of information transport, denote for Model A by $G_{\omega} \subset G$, the subgraph purely used to connect the information receiving component $\omega$ to the corresponding information sources for information transport. For example, in Figure~\ref{fig:explicit_information}, $\omega$ represents the component Objective Function and $G_{\omega}$ the blue colored subgraph. The implicit information transport results in a decrease of $\sum_{\omega} |V_{G_{\omega}}|$ nodes and $\sum_{\omega} |A_{G_{\omega}}|$ arcs in Model B in comparison to Model A.

For Model A, consider a generating unit $i \in I$ as the union of its $s$ subcomponents $i_1, \ldots i_s$. Since Model B does not represent the subcomponents, the number of nodes and arcs in Model B decreases by at least $s$, i.e., for each subcomponent its corresponding node and at least one connecting arc.

In total, Model B describes a subgraph of Model A with
\begin{equation}
    |V_A| \geq |V_B| + 2 |C|  + \sum_{\omega} |V_{G_{\omega}}| + s
\end{equation}
and 
\begin{equation}
    |A_A| \geq |A_B| + 2  |C|  + \sum_{\omega} |A_{G_{\omega}}| + s.
\end{equation}

\subsection{Comparison of MIP size}\label{section:comparison_mip}

Let $n_A, n_B \in \mathbb{N}$ the number of variables and $m_A, m_B \in \mathbb{N}$ the number of constraints in Model A and Model B.
Let $X_A, X_B$ be the set of flow variables in Model A and Model B, respectively.

Regarding the modeling of the flow variables $X$ either with a focus on nodes (Model A) or with a focus on arcs (Model B), the difference in the number of variables is given by $|X_A| = 2 \cdot |X_B|$. Figure~\ref{fig:flow_vars} depicts, how two variables in Model A can be represented as one variable in Model B.
Furthermore, for each arc $(v,w) \in A_A$, an identity constraint $x^{out}_v = x^{in}_w$ is necessary for the flow conservation. Therefore, the difference in the number of constraints is given by $m_A = m_B + |A_A|$.

Depending on the specific implementation, using ports on nodes, as in Model A, can result in unnecessary generated variables and constraints, even if the corresponding ports are not operated. Since the implementation can mitigate this issue by restricting the generation of variables and constraints by conditional queries to the active ports, we do not count for these extra variables and constraints.

The union of several conversion processes in Model B can decrease the size of the MIP significantly when several piece-wise linear functions are combined since modeling those processes usually relies on many variables and constraints. Denote the number of eliminated piecewise linear functions by $\mu$ and the number of reference points per piecewise linear function by $\rho$. This results in a decrease of $2\mu(\rho - 1)$ variables and $\mu(2\rho - 1)$ constraints.

In total, the difference in the size of the MIP in Model A and B is given by
\begin{equation}
    n_A = n_B + 2 \cdot |A_A| + 2\mu(\rho - 1) \geq n_B + 2 \cdot (|A_B| + 2 \cdot |C|  + \sum_{\omega} |A_{G_{\omega}}| + s) + 2\mu(\rho - 1)
\end{equation}
and
\begin{equation}
    m_A = m_B + |A_A| + \mu(2\rho - 1) \geq n_B + |A_B| + 2 \cdot |C|  + \sum_{\omega} |A_{G_{\omega}}| + s + \mu(2\rho - 1).
\end{equation}

\subsection{Mathematically equivalence}

Most of the considered decisions yield identical solutions for both models. For the proof of this identity, consider the notion of edge contractions in graphs.

\begin{definition} \label{def} 
    \parencite{Diestel2025}
    Let $e = uv$ be an edge of a graph $G = (V,E)$. Then $G/e$ is the graph obtained from G by \emph{contracting} the edge $e$ into a vertex $v_e$, i.e. $G/e = (V', E')$ with $V' = (V \setminus \{u,v\})\cup \{v_e\}$ and $E' = \{wz \in E \mid \{w,z\} \cap \{u,v\} = \emptyset\} \cup \{v_ez \mid uz \in E \setminus \{e \} \text{ or } vz \in E \setminus \{e\}\}.$
\end{definition} 

Definition~\ref{def} describes the notion of contraction to a general undirected graph. However, it can easily be adapted to our setting of a directed graph. 
We first see, that contracting the in- and outgoing arcs of balance nodes with degree two in the underlying network graph does not change the solution set.

\begin{lemma}
    Let $G = (V,A)$ be a network graph and let $u, w \in V$ and $v \in B \subset V$ a balance node such that $(u,v),(v,w) \in A$ are the only ingoing and outgoing arcs of $v$.  Contract the arcs $(u,v)$ and $(v,w)$ such that $G' = (V, (A \cup \{(u,w) \} ) \setminus \{ (u,v), (v,w) \}$. Then $x_G$ is an optimal solutions to $MIP_G$ (MIP based on graph $G$) if and only if $x_{G'}$ ($x_G$ restricted to $G'$) is an optimal solution to $MIP_{G'}$.
\end{lemma}

\begin{proof}
     Since the solution to the unit commitment problem is uniquely determined by its flow variables, we restrict the proof to $X$. The only constraint involving $v$ is $x_{in,v} = x_{out,v}$. Due to flow conservation on arcs, we know that $x_{in, v} = x_{out, u}$ and $x_{out, v} = x_{in, w}$, hence $x_{out, u} = x_{in, w}$. Since the contraction do not change the values of $x_{out, u} $ and $x_{in, w}$, the rest of the solution is not affected by the contraction.
\end{proof}

The iterative application of this result leads to an equivalence of Model A and B regarding container structures and additional transport of information.
Furthermore, merging sequential generation units in the underlying network graph does not change the solution set.

\begin{lemma}
    Let $G = (V,A)$ be a network graph and let $u, z \in V$ and $v, w \in I$ such that $(u,v),(v,w)$, $ (w,z) \in A$.  Let $G'$ be the result of the contraction of the arcs $(u,v), (v,w)$ and $(w,z)$ in $G$. Then $x_G$ is an optimal solution to $MIP_G$ if and only if $x_{G'}$ ($x_G$ restricted to $G'$) is an optimal solution to $MIP_{G'}$.
\end{lemma}

\begin{proof}
    We restrict the proof to the case, were $(u,v),(v,w), (w,z) \in A$ are the only ingoing and outgoing arcs of $v$ and $w$, i.e $G' = (V_G, (A_G \cup \{(u,vw), (vw,z) \} ) \setminus \{ (u,v), (v,w), (w,z) \}$. The more general result follows analogously. Since the solution to the unit commitment problem is uniquely determined by its flow variables, we restrict the proof to $X$. The flow conversion constraints involving $v$ and $w$ are $\varphi_1(x_{in,v}) = x_{out,v}$ and $\varphi_2(x_{in,w}) = x_{out,w}$. Due to flow conservation on arcs, we know that $x_{out, v} = x_{in, w}$ and hence $\varphi_2(\varphi_1(x_{in,v})) = \varphi_2 \circ \varphi_1(x_{in,v}) = x_{out,w}$. Possible ramping constraints and starting costs of the unit can be merged by eliminating the weaker constraint, respectively.
\end{proof}

This result shows the equivalence of Model A and B regarding the subcomponent structure. We further see that the modeling of flow variables in Model A and B is equivalent with the following result.

\begin{lemma}
    Given a graph $G = (V,A)$ and a corresponding network flow problem. Then the following is equivalent: 
\begin{enumerate}
    \item For all $v \in V$ define for each incoming arc a variable $x_{in,v}\in \mathbb{R}_{\geq 0}$ and each outgoing arc a variable $x_{out,v}\in \mathbb{R}_{\geq 0}$, representing the incoming and outgoing flow of $v$, respectively. For each arc $(v,w) \in A$ add the flow conservation constraint $x_{out, v} = x_{in, w}$.
    \item For each arc $(v,w) \in A$, define a variable $x_{(v,w)} \in \mathbb{R}_{\geq 0}$ representing the value of flow on arc $(v,w)$.
\end{enumerate}
\end{lemma}

\begin{proof}
    Let $(v,w) \in A$. Since $x_{out, v} = x_{in, w}$, a representation by $x_{(v,w)}$ (or vice versa) with $x_{(v,w)} = x_{out, v} = x_{in, w}$ is feasible.
\end{proof}

The specification of node types and the usage of ports for connections between nodes only influence the ease of implementation and handling and have, therefore, no influence on the equivalence of the models. The equivalence of detail of technical constraints for specific technologies highly depends on the specific constraints. If there is a need for mathematical equivalence, the reduction of constraints should exclusively be restricted to redundant constraints. In total, we have seen that the dimensions considered in Model A and Model B do not influence the solutions to optimization problems.

\section{Computational study} \label{section:study}

The following computational study aims to depict the difference in the problem size and the solution progress for Model A and Model B in a real-world example.

\subsection{Instances}

The instances are based on a realistic subnetwork of the Berlin district heating system. The subnetwork is characterized by a centralized structure with four heat generating units located on three sites.
The primary heat generation is given by two plants: One is a heating plant operated by gas, with a share of $81.6$~\% of the total capacity. The other is a combined heat and power plant (CHP) operated by biomass, with a share of $16.8$~\% of the total capacity. 
Those primary heat generation units are supplemented by two smaller CHPs: One is operated by gas, with a share of $1.3$~\%, and the other is operated by biogas, with a share of $0.2$~\% of the total heat generation capacity. 
A storage facility can store up to $84,74~\%$ of the maximal heat that could be generated per hour. 
We consider a time horizon of one month with a 4-hour time step. 
The instances involve three hierarchical objectives: the minimization of costs ($f_1$), the minimization of CO$_2$ emissions ($f_2$) and the maximization of heat ouput of CHPs ($f_3$). This study involves three scenarios representing one month in the onset, midseason, and conclusion heating seasons.

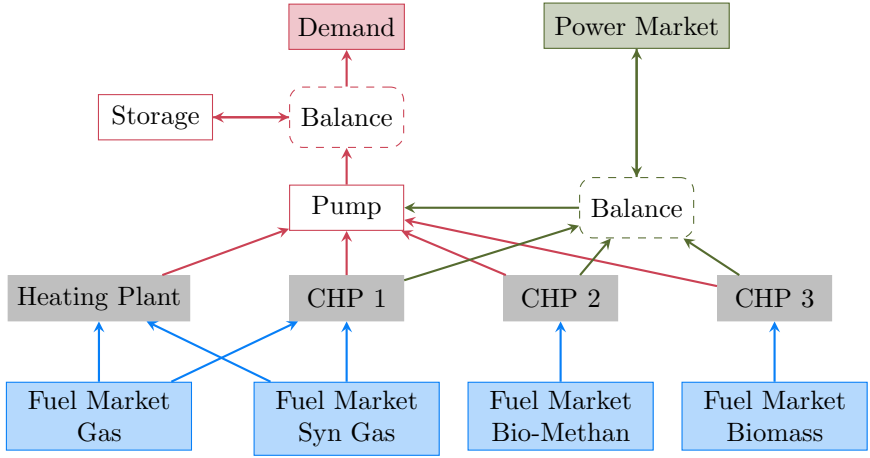
\begin{figure}[h!]
    \centering
    \tikzstyle{container} = [rectangle, 
minimum width=1.5cm, 
minimum height=0.6cm, 
text centered, 
draw=myred,
fill=myred!30]
\tikzstyle{containergreen} = [rectangle, 
minimum width=1.5cm, 
minimum height=0.6cm, 
text width=2.2cm,
text centered, 
draw=mygreen,
fill=mygreen!30]
\tikzstyle{container1} = [rectangle, 
minimum width=1.5cm, 
minimum height=0.6cm, 
text centered, 
draw=lightgray,
fill=lightgray]
\tikzstyle{container2} = [rectangle, 
minimum width=1.5cm, 
minimum height=0.6cm, 
text width=2.2cm,
text centered,
draw=myblue,
fill=myblue!30]
\tikzstyle{container3} = [rectangle, 
minimum width=1.5cm, 
minimum height=0.6cm, 
text centered, 
draw=myred]
\tikzstyle{largecontainer} = [rectangle, 
minimum width=10cm, 
minimum height=4.9cm, 
text centered, 
draw=myblue]
\tikzstyle{oval} = [rectangle, 
minimum width=1.5cm, 
minimum height=0.8cm, 
text centered,
dashed,
draw=myred,
rounded corners=5pt]
\tikzstyle{ovalgreen} = [rectangle, 
minimum width=1.5cm, 
minimum height=0.8cm, 
text centered,
dashed,
draw=mygreen,
rounded corners=5pt]
\tikzstyle{enhanced} = [circle, 
minimum width=1cm, 
minimum height=1cm, 
text centered, 
dashed,
draw=myblue]
\tikzstyle{dummy} = [circle, 
text centered]
\tikzstyle{dummy2} = [rectangle, 
text centered]
\tikzstyle{arrow} = [thick,->,>=stealth]

    \begin{tikzpicture}[node distance=1.2 cm]

    \node (demand) [container] {Demand};
    \node (balance) [oval, below of= demand] {Balance};
    \node (storage) [container3, left=1.0cm of balance] {Storage}; 
    \node (pumpe) [container3, below of=balance] {Pump}; 
    \node (H2) [container1, below of=pumpe] {CHP 1}; 
    \node (H1) [container1, left=1.3cm of H2] {Heating Plant}; 
    \node (KWK1) [container1, right=1.3cm of H2] {CHP 2}; 
    \node (KWK2) [container1, right=1.3cm of KWK1] {CHP 3};
    \node (fuel1) [container2, below=0.8cm of H1] {Fuel Market \\ Gas};
    \node (fuel2) [container2, below=0.8cm of H2] {Fuel Market \\ Syn Gas};
    \node (fuel3) [container2, below=0.8cm of KWK1] {Fuel Market \\ Bio-Methan};
    \node (fuel4) [container2, below=0.8cm of KWK2] {Fuel Market \\ Biomass};
    \node (balance2) [ovalgreen, right=2.3cm of pumpe] {Balance};
    \node (power) [containergreen, above=1.7cm of balance2] {Power Market};

    \draw [arrow, myred] (balance) --  (demand);
    \draw [arrow, myred] (balance) --  (storage);
    \draw [arrow, myred] (storage) --  (balance);
    \draw [arrow, myred] (pumpe) --  (balance);
    \draw [arrow, myred] (H1) --  (pumpe);
    \draw [arrow, myred] (H2) --  (pumpe);
    \draw [arrow, myred] (KWK1) --  (pumpe);
    \draw [arrow, myred] (KWK2) --  (pumpe);
    \draw [arrow, myblue] (fuel1) --  (H1);
    \draw [arrow, myblue] (fuel1) --  (H2);
    \draw [arrow, myblue] (fuel2) --  (H1);
    \draw [arrow, myblue] (fuel2) --  (H2);
    \draw [arrow, myblue] (fuel3) --  (KWK1);
    \draw [arrow, myblue] (fuel4) --  (KWK2);

    \draw [arrow, mygreen] (KWK1) --  (balance2);
    \draw [arrow, mygreen] (KWK2) --  (balance2);
    \draw [arrow, mygreen] (H2) --  (balance2);
    \draw [arrow, mygreen] (balance2) --  (pumpe);
    \draw [arrow, mygreen] (balance2) --  (power);
    \draw [arrow, mygreen] (power) --  (balance2);

    \end{tikzpicture}
    \caption{Simplified illustration of the instance of a realistic subnetwork of the Berlin district heating system (fuel (blue), heat (red), power (green)). }
    \label{fig:MV}  
\end{figure}

\subsection{Computational results}

We solve the multi-objective instances with lexicographic optimization to obtain comparable solutions for both models, i.e. for the three objective functions $f_1, f_2, f_3$ compute iteratively:
\begin{enumerate}
    \item $f_1^* = \min f_1$
    \item $f_2^* = \min f_2$ s.t. $f_1 \leq f_1^*$
    \item $f_3^* = \min f_3$ s.t. $f_1 \leq f_1^*$ and $f_2 \leq f_2^*$
\end{enumerate}

All computation where performed using a single Intel(R) Xeon(R) CPU E3-1245 v6 @ 3.70GHz. The commercial solver Gurobi Optimizer version 11.0.0 \parencite{gurobi} 
was applied to each single optimization problem. 

\begin{table} [h!]
    \centering
    \resizebox{\linewidth}{!}{
    \begin{tabular}{llrrrrrr}
    \hline
        & &        & Model A &   &   & Model B & \\
        & & onset       & midseason & conclusion  & onset  & midseason & conclusion\\
    \hline
    \hline
      Graph &\# nodes & & 128& & & 33& \\
     &\# arcs & & 171& & & 38& \\
    \hline
    \hline
     MIP &\# variables total & 74029 & 74029 & 71641& 16740 & 16740 & 16200\\
     &\# constraints & 70309 & 70309 & 68041 & 17102 & 17102 & 16550\\
    \hline
    Lex. 1 &\# variables after presolve & 3528 & 3510 &3417 &3964 & 4082& 3901\\
    &\# constraints after presolve & 4262 & 4229 &4127 & 2837 & 2967 & 2796\\
    \hline
    Lex. 2 &\# variables after presolve & 4639 & 4718 &4313 & 4979& 4964& 4813 \\
    &\# constraints after presolve & 4076 & 4415 & 3768& 4398 & 4376 & 4248 \\
    \hline
    Lex. 3 &\# variables after presolve & 4643 & 4746  & 4493& 4979& 4964& 4813\\
    &\# constraints after presolve & 4081 & 4184 &3949 & 4399 & 4377 & 4249\\
    \hline
    Time & total presolve time [s] & 0.23  & 0.20  & 0.31 & 0.16 &  0.21 & 0.16 \\
    & total runtime [s] & 0.72  & 0.47  & 1.27 & 0.67 &  0.37& 0.93 \\
    & presolve disabled [s] & 8.72 & 3.48 & 30.07 & 2.31 & 1.04& 3.55\\
    \hline
    \end{tabular}}
    \caption{Size of underlying network structure and the initial MIP. Presolve results depict the size of MIP in each lexicographic optimization call. The runtimes are measured as the average of five runs.}
    \label{tab:table2}
\end{table}

Table~\ref{tab:table2} shows the results for the three instances of Model A and Model B.
Model B shows a significantly smaller graph size, with the number of nodes decreased by $74.2\%$ and the number of arcs decreased by $77.8\%$. For all instances, the initial mixed integer program also decreases with Model B in the number of variables by $77.4\%$ and in the number of constraints by $75.7\%$. Those results confirm the advantage of Model B regarding a less complex graph and the size of the resulting MIP. 
However, the average MIP size for each lexicographic optimization call after the respective presolve procedure shows a similar size for both models, even often smaller for Modal A. This indicates that the presolve method of Gurobi will automatically and successfully deal with the redundant formulations of Model A. 

The runtimes for Model B show an average decrease of $19.9\%$ in comparison to Model A. The longer runtimes for Model A can largely be attributed to the additional effort in the presolving stage. The presolving time in Model B is on average $28.4\%$ shorter than in Model A.
If the presolve is disabled, the runtime of model B even shows an average decrease of $83.6\%$

\section{Conclusions}
\label{sec:conclusions}

For multi-energy system optimization frameworks based on mixed integer programs, we described a set of modeling dimensions that are usually not discussed in their classification.
We defined two models based on those dimensions with diametrical characteristics, one focusing on human readability and the other on mathematical flexibility. The first model showed a strong focus on the individual units of the network, which is a substantial advantage for visual frameworks and the inclusion of technology-specific modeling. The second model represents a strong abstraction of the structures and an arc-based flow model, as often viewed in graph-theoretical literature, and is, therefore, beneficial for the application of algorithms. It further comes with a substantial decrease in the size of the underlying graph and the resulting MIP. A study on a realistic instance of an urban district heating network confirmed the advantage in size. However, the applied commercial solver recognized the redundancies in the presolve step. Therefore, the advantage in problem size only impacts the presolve time and storage use. All of the considered dimensions are highly dependent on what audience and which use the framework is aimed at. 
Future research could include the investigation of further overlooked dimensions, a rigorous survey of existing frameworks about the considered and further dimensions, and a study on a more extensive set of instances. 
The presented work aims to increase the awareness of multi-energy system optimization framework developers about the impact of specific design decisions.

\section*{Nomenclature}
\label{sec:symbole_abk}

\subsection*{Abbreviations}
\begin{tabular}{@{}p{3cm}l}
    MIP & Mixed-integer program\\
\end{tabular}

\subsection*{Superscripts and Subscripts}
\begin{tabular}{@{}p{3cm}l}
    $r$ & Resource \\
	$i/k$ & Generating/ storage unit \\
    $in/ out$ & Ingoing/ outgoing flow \\
    $t$ & Time step \\
\end{tabular}

\subsection*{Greek Symbols}
\begin{tabular}{@{}p{3cm}l}
    $\varphi, \psi$ & Resource conversion map \\
\end{tabular}

\subsection*{Latin Symbols}
\begin{tabular}{@{}p{3cm}l}
	$G = (V,A)$ & Graph $G$ with set of vertices $V$ and set of arcs $A$\\
    $u,v,w$ & Nodes in graph $G$ \\
    $e = vw,  a = (v,w)$ & Edges and arcs in graph $G$ \\
    $R$ & Set of resources \\
    $I, K$ & Set of generating units, set of storage units\\
    $B$ & Set of balance nodes\\
    $T$ & Set of time steps\\
    $C$ & Set of container structures\\
    $x, z, s, h$ & Variables for resource flow, operation, status and storage level\\
    $p, e$ & Variables for purchased and sold resources \\
    $d, a, c$ & Parameter vector for demand, coefficients and capacities \\
    $f_1, f_2, f_3$ & Ordered objective functions\\
    $f^*$ & optimal value of objective function $f$ \\
    $n, m$ & Number of variables and constraints\\
    $s$ & Number of subcomponents\\
    $x_i, \overline{y}_i, x^*_i$ & Source series for characteristic curve, $i = 0, \dots, n$\\
    $y_i, z_i, z^*_i$ & Target series for characteristic curve, $i = 0, \dots, n$ \\
\end{tabular}

\printbibliography

\section*{Acknowledgement}
The work for this article has been conducted in the Research Campus MODAL funded by the German Federal Ministry of Education and Research (BMBF) (fund numbers 05M14ZAM, 05M20ZBM). The work has been supported by the German National Science Foundation (DFG) Cluster of Excellence MATH+.

\end{document}